\documentclass[runningheads]{llncs}
\usepackage{amsmath, amssymb, amsfonts}
\usepackage{xcolor}
\usepackage{booktabs}

\usepackage{listings}
\usepackage{xspace}

\usepackage[ngerman,british]{babel}

\usepackage[style=alphabetic,
   sorting=nty,
   backend=biber,
   maxbibnames=99,
   doi=false,
   isbn=false,
   backref=false]{biblatex}
\usepackage[tracking, kerning, spacing]{microtype}
\microtypecontext{spacing=nonfrench}

\usepackage{tikz}
\usepackage{tikz-3dplot}
\usetikzlibrary{matrix,fit,positioning,shapes.geometric,patterns,backgrounds}
\usetikzlibrary{decorations.pathreplacing}
\usetikzlibrary{arrows,calc}
\tikzstyle{bigbox} = [draw=blue!50, thick, rounded corners, rectangle]
\tikzset{
>=stealth'
}

\usepackage{algorithm}
\usepackage{algpseudocode}

\definecolor{larscolor}{rgb}{.7,.5,0}

\definecolor{martacolor}{rgb}{.7,0,.5}

\usepackage[normalem]{ulem}

\newcommand{\polymake}{{\texttt{polymake}}\xspace}

\newcommand{\zram}{\textbf{ZRAM}\xspace}

\usepackage{aliascnt}

\usepackage{hyperref}

\newcounter{internal}[section]
\renewcommand{\theinternal}{\thesection.\arabic{internal}}

\newaliascnt{intcor}{internal} 
\newaliascnt{intconj}{internal} 
\newaliascnt{intlemma}{internal} 
\newaliascnt{intdef}{internal} 
\newaliascnt{intex}{internal} 
\newaliascnt{intprop}{internal} 
\newaliascnt{intrem}{internal} 
\newaliascnt{intthm}{internal}

\newcommand*{\dupcntr}[2]{%
   \expandafter\let\csname c@#1\expandafter\endcsname\csname c@#2\endcsname
}
\dupcntr{figure}{internal}
\renewcommand{\thefigure}{\theinternal}

\newtheorem{prop}[intprop]{Proposition}

\definecolor{mcolor}{rgb}{.4,0,.4}
\newcommand{\macrocolor}[1]{\textcolor{black}{#1}}

\newcommand{\RR}{\macrocolor{\mathbb R}}

\newcommand{\scalp}[1]{\macrocolor{\langle} #1 \macrocolor{\rangle}}
\newcommand{\supportcone}{{\macrocolor{\mathcal{S}}}}
\newcommand{\subdivision}{{\macrocolor{\Sigma}}}
\newcommand{\zonotope}{{\macrocolor{\mathcal{Z}}}}
\DeclareMathOperator{\pc}{\macrocolor{\mathcal{PC}}}
\DeclareMathOperator{\cone}{\macrocolor{cone}}
\DeclareMathOperator{\maxcones}{\macrocolor{maxcones}}
\DeclareMathOperator{\signature}{\macrocolor{sig}}

\begin{document}

\title{Hyperplane arrangements in \polymake}
\author{
Lars Kastner \and
Marta Panizzut
}
\authorrunning{L. Kastner and M. Panizzut}
\institute{Technische Universit\"at Berlin, Chair of Discrete Mathematics/Geometry,
Stra\ss e des 17. Juni 136, 
10623 Berlin 
\email{\{kastner,panizzut\}@math.tu-berlin.de}}
\maketitle              %
\begin{abstract}
Hyperplane arrangements form the latest addition to the zoo of combinatorial objects dealt with by \polymake.  We report on their implementation 
and on a algorithm to compute the associated cell decomposition.  The implemented algorithm performs significantly  better than brute force alternatives, as it  requires less convex hulls computations.

\end{abstract}

\section{Introduction}

Hyperplane arrangements are ubiquitous objects appearing in different areas of
mathematics such as discrete geometry, algebraic combinatorics and algebraic
geometry. A common theme is to understand the combinatorics and the topology of
the cells in the complement of the arrangement. Combinatorics and its
connections to other areas of mathematics are the focus of the software
framework \polymake \cite{DMV:polymake}, hence hyperplane arrangements form an almost mandatory
addition to the objects available. We will discuss the implementation, such as
the datatypes and properties, as well as some basic algorithms for analyzing
hyperplane arrangements.

One of the main advantages of \polymake are its various interfaces to other
software. This allows keeping the codebase slim, while using powerful software
that developed by experts from other fields. Still \polymake provides basic
algorithms for many tasks, in case other software is not available. Hence the
idea of the hyperplane arrangements is to provide a datatype with basic
functionality as basis for future interfaces to other software, e.g. to \zram
\cite{zram} for computing the cell decomposition from the hyperplanes.
Nevertheless, the \polymake implementation of hyperplane arrangements comes
with a basic algorithm for computing the associated cell decomposition that
performs significantly better than brute force alternatives. Thus, we will
discuss the main ideas of this algorithm in this article as well.

The combinatorics of hyperplane arrangements in real space is linked to zonotopes. 
Each arrangement endows the support space with a fan structure which is the  normal fan of a zonotope. Each hyperplane subdivides the space in two halfspaces. Therefore we can encode relative positions of points with respect to the arrangement. In other words, hyperplane arrangements are examples of (oriented) matroids.  Moreover, the hyperplanes in an arrangement can be seen as  mirrors hyperplanes of a reflection group.  
 
An interesting application is in Geometric Invariant Theory. GIT constructs quotients of algebraic varieties modulo group actions. The quotients depend on the choice of a linearized ample line bundle. Variation of geometric invariant theory quotients studies how quotients vary when changing the line bundle. Under some hypothesis the classes of equivalent quotients are convex subsets, called chambers. The walls among chambers are defined by certain hyperplane arrangements, see \cite[Example 3.3.24]{DolgachevHu}.

\section{Main Definitions}
We begin with the basic definitions in the theory of hyperplane arrangements following our implementation in \polymake.

\begin{definition}
A \emph{hyperplane arrangement} $H=(H,\supportcone_H)$  in $\mathbb{R}^d$ is given by the following data:
\begin{enumerate}
\item A finite set of linear forms encoding hyperplanes $H = \big\{h \in \RR^d \setminus \{0\}\big\}$ and
\item A polyhedral cone $\supportcone_H\subseteq\RR^d$ which we call the \emph{support
cone}.
\end{enumerate}
Given a hyperplane arrangement $H$, the \emph{induced fan} $\subdivision_H$ is a fan
with support $\supportcone_H$ given by subdividing $\supportcone_H$ along all $ \{x \in \RR^d \, | \, \scalp{x,h} =0\}$ for $h \in H$. 

Every hyperplane in the arrangement subdivides the space into two halfspaces
\[
h^+\ :=\ \{x\in\RR^d \, | \, \scalp{x,h} > 0\}
\mbox{ and }
h^-\ :=\ \{x\in\RR^d\, | \, \scalp{x,h} < 0\}.
\]
\end{definition}

We remark that in the definition  we allow duplicate hyperplanes, but from each hyperplane arrangement we can construct a reduced one. 
Let $H$ be a hyperplane arrangement given by the hyperplanes $\{h_1, h_2, \dots, h_n\}$. The \emph{reduced hyperplane arrangement} $H_{\text{red}}$ has the same support cone as $H$ and  $h_i \in H_{\text{red}}$ if and only if $h_i \not = \lambda b$ , for any $\lambda \in \RR$ and any $b \in \{h_1,  \dots, h_{i-1}\}$.

To a hyperplane arrangement $H=\{h_1,\ldots,h_n\}\subseteq\RR^d$ we associate
the polytope
\[
\zonotope_H\ :=\ \sum_{i=0}^n \, [-h_i,h_i]+\supportcone_H^{\vee},
\]
the Minkowski sum of all the line segments $[-h_i,h_i]$ and the dual support
cone $\supportcone_H^{\vee}$.  If $\supportcone_H=\RR^d$, then $\supportcone_H^{\vee}=0$ and $\zonotope_H$ is a zonotope.

\begin{remark}
Often hyperplane arrangements are defined without a support cone, i.e. only for
the case $\supportcone_H=\RR^d$. The connection between intersecting
$\subdivision_H\cap\supportcone_H$ is done via taking the Minkowski sum
$\zonotope_H+\supportcone_H^{\vee}$ on the dual side. The main ingredient is
the fact that
\[
(\sigma+\tau)^\vee\ =\ \sigma^\vee\cap\tau^\vee
\]
holds for two cones $\sigma$ and $\tau$.
\end{remark}

\begin{prop}{{\cite[Thm. 7.16]{ziegler}}}
The fan $\subdivision_H$ is the normal fan of $\zonotope_H$.
\end{prop}

\begin{definition}
To a maximal cone of $\sigma\in\subdivision_H$ we associate its
\emph{signature}, which is a set 
$\signature(\sigma):=\big\{i\in\{1,\ldots,n\}\ |\ \sigma\subseteq \overline{h_i^-}\big\}$.
\end{definition}

\begin{example}\label{ex:2dim}
Let $H$ be given by
\[
H\ =\ \{(0,1),(1,1),(-2,1)\}\ \subseteq\RR^2.
\]
We will have a look at the induced fans for different support cones $\supportcone_H$. The
fan $\subdivision_H$ and the polytope $\zonotope_H$ are visualized in \autoref{fig:ex:2dim}
for varying $\supportcone_H$.

\begin{figure}
\def\picscale{.5}
\def\baseline{0cm}
\def\hypwidth{1pt}
\[
\begin{array}{r|ccc}
\supportcone_H & \RR^2 & \RR_{\ge 0}^2 & \cone\{(0,1),(1,-1)\}\\
\midrule
\subdivision_H &
\begin{tikzpicture}[scale=\picscale, baseline=\baseline]
\draw[color=black!40] (-2.3,-2.3) grid (2.3,2.3);
\fill[pattern color=black!20, pattern=north west lines] (-2.3,-2.3) rectangle (2.3,2.3); 
\draw[line width=\hypwidth] (-2.3/2,-2.3) -- (2.7/2,2.7);
\draw[line width=\hypwidth] (2.3,-2.3) -- (-2.7,2.7);
\draw[line width=\hypwidth] (-2.7,0) -- (2.3,0);
\node[above right] at (-2.3,2.3) {\tiny{$+$}};
\node[below left] at (-2.3,2.3) {\tiny{$-$}};
\node[above left] at (-2.3,0) {\tiny{$+$}};
\node[below left] at (-2.3,0) {\tiny{$-$}};
\node[above left] at (2.3/2,2.3) {\tiny{$+$}};
\node[above right] at (2.3/2,2.3) {\tiny{$-$}};
\node at (2,1) {\tiny{\textbf{$++-$}}};
\node at (0,-2) {\tiny{\textbf{$---$}}};
\end{tikzpicture}
&
\begin{tikzpicture}[scale=\picscale, baseline=\baseline]
\draw[color=black!40] (-2.3,-2.3) grid (2.3,2.3);
\fill[pattern color=black!40, pattern=north west lines] (0,0) rectangle (2.3,2.3); 
\draw[line width=\hypwidth] (-2.3/2,-2.3) -- (2.3/2,2.3);
\draw[line width=\hypwidth] (2.3,-2.3) -- (-2.3,2.3);
\draw[line width=\hypwidth] (-2.3,0) -- (2.3,0);
\end{tikzpicture}
&
\begin{tikzpicture}[scale=\picscale, baseline=\baseline]
\draw[color=black!40] (-2.3,-2.3) grid (2.3,2.3);
\fill[pattern color=black!40, pattern=north west lines] (2.3,-2.3) -- (0,0) -- (0,2.3) -- (2.3,2.3) -- cycle;
\draw[line width=\hypwidth] (-2.3/2,-2.3) -- (2.3/2,2.3);
\draw[line width=\hypwidth] (2.3,-2.3) -- (-2.3,2.3);
\draw[line width=\hypwidth] (-2.3,0) -- (2.3,0);
\end{tikzpicture}\\
\midrule
\zonotope_H &
\begin{tikzpicture}[scale=.4, baseline=\baseline]
\draw[color=black!40] (-3.3,-3.3) grid (3.3,3.3);
\fill (0,0) circle (4pt);
\fill[pattern color=black!40, pattern=north west lines] (1,-3) -- (3,-1) -- (3,1) -- (-1,3) -- (-3,1) -- (-3,-1) -- cycle;
\draw[line width=\hypwidth] (1,-3) -- (3,-1) -- (3,1) -- (-1,3) -- (-3,1) -- (-3,-1) -- cycle;
\end{tikzpicture}
&
\begin{tikzpicture}[scale=.4, baseline=\baseline]
\draw[color=black!40] (-3.3,-3.3) grid (3.3,3.3);
\fill (0,0) circle (4pt);
\fill[pattern color=black!40, pattern=north west lines] (-3,3.3) -- (-3,-1) -- (1,-3) -- (3.3,-3) -- (3.3,3.3) -- cycle;
\draw[line width=\hypwidth] (-3,3.3) -- (-3,-1) -- (1,-3) -- (3.3,-3);
\end{tikzpicture}
&
\begin{tikzpicture}[scale=.4, baseline=\baseline]
\draw[color=black!40] (-3.3,-3.3) grid (3.3,3.3);
\fill (0,0) circle (4pt);
\fill[pattern color=black!40, pattern=north west lines] (1,-3) -- (3.3,-3) -- (3.3,3.3) -- (-1+.3,3+.3) -- (-3,1) -- (-3,-1) -- cycle;
\draw[line width=\hypwidth] (-1+.3,3+.3) -- (-3,1) -- (-3,-1) -- (1,-3) -- (3.3,-3);
\end{tikzpicture}
\\
\midrule
\#\maxcones{\subdivision_H} & 6 & 2 & 3
\end{array}
\]
\caption{Visualization of $\subdivision_H$ and $\zonotope_H$ for \autoref{ex:2dim}}
\label{fig:ex:2dim}
\end{figure}
In each of the pictures, the support cone is indicated as the shaded area. The
structure of the fan $\subdivision_H$ depends heavily on the support cone $\supportcone_H$. In
particular, it is possible for hyperplanes to only intersect $\supportcone_H$ trivially
and thereby becoming irrelevant for $\subdivision_H$. Thus, one may loose
information when going from $H$ to $\subdivision_H$.

The labels at the hyperplanes in the first picture indicate which side
constitutes $h^+$, $h^-$ respectively. Using these one can read of the
signatures of the single cells, for example the cell $\sigma$ generated by the rays
$(1,0)$ and $(1,2)$ has signature $\signature(\sigma)=\{3\}$.
\end{example}

\begin{remark}
Reduced hyperplane arrangements are examples of oriented matroids. The ground set is the collection $H$ of hyperplanes and the signatures are the covectors. 
\end{remark}

\subsection{Affine hyperplane arrangements}

An affine hyperplane arrangement is usually given by a finite set of affine
hyperplanes:
\[
H_{\text{aff}}\ :=\ \{[a,b]\in\RR^{d}\times \RR\}.
\]
The whole space $\RR^d$ is then subdivided along the hyperplanes 
\[
\{x\in\RR^d\, |\, \scalp{a,x}=b\}, \mbox{ for all }[a,b]\in H,
\]
resulting in a polyhedral complex $\pc_{H_{\text{aff}}}\subseteq\RR^d$.

Analogously to the connection between polytopes and cones, or polyhedral
complexes and fans, every affine hyperplane arrangement gives rise to a
(projective) hyperplane arrangement by embedding it at height $1$:
\[
H_{\text{proj}}\ :=\ \{[-b,a]\ |\ [a,b]\in H\}.
\]
If we intersect the fan $\subdivision_{H_{\text{proj}}}$ with the affine hyperplane
$[x_0=1]\subseteq\RR^{d+1}$, the resulting polyhedral complex is isomorphic to
$\pc_{H_{\textrm{aff}}}$, via the embedding $\RR^d\to\RR^{d+1}$, $x\mapsto
[1,x]$.

The support cone allows one to deal with affine hyperplanes computationally. Set
\[
\supportcone_{H_{\text{proj}}}\ :=\ \{[x_0,x_1,\ldots,x_{d}]\in\RR^{d+1}\ |\ x_0\ge 0\},
\]
then the maximal cones of $\subdivision_{H_{\text{proj}}}$ are in one-to-one correspondence
with the maximal cells of $\pc_{H_{\text{aff}}}$. In particular, \polymake can
interpret $\subdivision_{H_{\text{proj}}}$ as a polyhedral complex via the embedding
mentioned above, and this polyhedral complex will be exactly $\pc_{H_{\text{aff}}}$.

\begin{example}
As a simple example, choose the following hyperplanes in $\RR^1$:
\[
\begin{array}{ccc}
x_1=-1, & x_1=0, & x_1=2.
\end{array}
\]
The associated hyperplanes of $H_{\text{proj}}$ in $\RR^2$ are exactly those of the
hyperplane arrangement from \autoref{ex:2dim}. For $\supportcone_H$ we choose the cone
$\{x_0\ge 0\}$, then $H_{\text{aff}}$ will be at height one.
\[
\begin{tikzpicture}[scale=.7]
\draw[color=black!40] (-2.3,-2.3) grid (2.3,2.3);
\fill[pattern color=black!40, pattern=north west lines] (0,-2.3) rectangle (2.3,2.3);
\draw[line width=1pt] (-2.3/2,-2.3) -- (2.3/2,2.3);
\draw[line width=1pt] (2.3,-2.3) -- (-2.3,2.3);
\draw[line width=1pt] (-2.3,0) -- (2.3,0);
\draw[line width=2pt] (1,-2.3) node[below] {$\pc_{H_{\text{aff}}}$} -- (1,2.3);
\fill (1,2) circle(4pt);
\fill (1,0) circle(4pt);
\fill (1,-1) circle(4pt);
\end{tikzpicture}
\]
The induced affine hyperplane arrangement is indicated by the dots and thick line. It is one dimensional and the associated polyhedral complex $\pc_{H_{\text{aff}}}$ has four maximal cells.
\end{example}

\begin{example}
\begin{lstlisting}
fan > $HA = new HyperplaneArrangement(HYPERPLANES=>[[0,1],[1,1],[-2,1]],"SUPPORT.INEQUALITIES"=>[[1,0]]);

fan > $HA->CELL_DECOMPOSITION->RAYS; # Force computation

fan > $pc = new PolyhedralComplex($HA->CELL_DECOMPOSITION);

fan > print "(".join("),(",@{$pc->VERTICES}).")\n";
(0 -1),(0 1),(1 -1),(1 0),(1 2)

fan > print join(",",@{$pc->MAXIMAL_POLYTOPES})."\n";
{0 2},{1 4},{2 3},{3 4}
\end{lstlisting}
\end{example}

\section{Implementation}

Hyperplane arrangements are implemented in the software \polymake as a new object
\texttt{HyperplaneArrangement}, which is derived from the already existing object 
 \texttt{VectorConfiguration}. Besides the existing properties of
\texttt{VectorConfiguration} it has been augmented with the following
properties and methods.
\begin{enumerate}
\item \textbf{\texttt{HYPERPLANES}} 
   A matrix containing the hyperplanes as rows, this is just an override of the
   property \texttt{VECTORS} of \texttt{VectorConfiguration}
\item \textbf{\texttt{SUPPORT}} 
   A \polymake \texttt{Cone}, denoting the support $\supportcone_H$.
\item \textbf{\texttt{CELL\_DECOMPOSITION}} 
   A \polymake \texttt{PolyhedralFan}, the cell decomposition~$\subdivision_H$.
\item \textbf{\texttt{CELL\_SIGNATURES}}
   A \texttt{Array<Set<Int>>}, the $i$-th set in the array contains the indices of
   hyperplanes evaluating negatively on the $i$-th maximal cone of
   \texttt{CELL\_DECOMPOSITION}.
\item \textbf{\texttt{signature\_to\_cell}}
   Given a signature as \texttt{Set<Int>}, get the maximal cone with this
   signature, if it exists.
\item \textbf{\texttt{cell\_to\_signature}}
   Given a cell, a maximal cone of \texttt{CELL\_DECOMPOSITION}, determine its
   signature.
\end{enumerate}

\subsection{Cell decomposition algorithm}
Given $H=\{h_1,\ldots,h_n\}$, we want to compute the subdivision of $\supportcone_H$
induced by the hyperplanes, the induced fan $\subdivision_H$. This means, we want to
find all the rays and maximal cones of $\subdivision_H$. In terms of the zonotope
$\zonotope_H$, this is equivalent to knowing the facets and vertices of $\zonotope_H$, see \cite{Fukuda, GritzmannSturmfels}.  The
facet directions of $\zonotope_H$ are the rays of $\subdivision_H$. For very vertex of $\zonotope_H$
we get a maximal cone by determining which facets contain it.

The brute force approach is to loop over all possible signatures in

$s\in 2^{0,\ldots,n}$ and for every signature $s$ to build the cone
\[
\bigcap_{i\in s}\overline{h_i^-}\ \cap\ 
\bigcap_{i\notin s}\overline{h_i^+}\ \cap\ 
\supportcone_H.
\]
For comparing the different algorithms, we count the number of times they have
to perform a convex hull computation for converting a signature to a cone.
There are $2^{n}$ signatures, so we have to perform $2^{n})$ convex hull
computations.  As we saw in \autoref{ex:2dim}, it can happen that some
hyperplanes are irrelevant, either completely or just for single cells.
Furthermore, in \autoref{ex:2dim} the maximum number of two-dimensional cones
we got was six, however we would have to compute eight intersections with the
brute force approach regardless.

\begin{remark}
This brute force approach is in some ways parallel to the brute force approach
for computing the Minkowski sum making up $\zonotope_H$, by taking considering
all possible sums of the endpoints of the line segments.  One arrives at
$2^{n}$ points whose convex hull is $\zonotope_H$. There are several ways to
go on: Either attempt a massive convex hull computation directly, or check each
point individually whether it is a vertex.
\end{remark}

Our approach is to first find a full-dimensional cone $\sigma$ of $\subdivision_H$
and then to flip hyperplanes in order to compute its neighbors. First take a
facet $f$ of $\sigma$, then set
\[
\signature'\ :=\ (\signature(\sigma)\setminus\{i\in\signature(\sigma)\ |\ h_i||f\})\cup\{i\notin\signature(\sigma)\ |\ h_i||f\})
\]
This is the signature of the cell neighboring $\sigma$ at facet $f$, so we can
use it to determine the rays of the neighboring cell. Finding the neighbors of
a cell allows one to traverse the dual graph of the fan $\subdivision_H$.
Taking the support cone $\supportcone_H$ into account just requires some minor
tweaks, like ignoring facets of $\sigma$ that are also facets of
$\supportcone_H$. By storing signatures one can avoid recomputation of cones.

To find a starting cone, one selects a generic point $x$ from $\supportcone_H$. A generic
point will be contained in a maximal cone, this maximal cone will be
\[
\sigma(x)\ :=\ \bigcap_{i\ |\ x\in h_i^+}\overline{h_i^+}\ \cap\ \bigcap_{i\ |\ x\in h_i^-}\overline{h_i^-}.
\]
The point $x$ may be contained in some hyperplanes, but these hyperplanes are
exactly those that also contain the entire $\supportcone_H$. Using this
approach we would do one convex hull computation per maximal cone, arriving at
$\#\maxcones(\subdivision_H)$ convex hull computations.

\begin{remark}
As the fan $\subdivision_H$ is polytopal, there is a reverse search structure
on it, corresponding to the edge graph of the zonotope $\zonotope_H$. This has
already been exploited by Sleumer in \cite{sleumer99} using the software
framework \cite{zram}. Reverse search allows for different kinds of
parallelisation and it would be interesting to study the performance of
budgeted reverse search \cite{mplrs, mts_tutorial} on this particular problem.
Note that the dual problem, finding the vertices of $\zonotope_H$, is equally
hard, as it is a Minkowski sum with potentially many summands. We refer the
reader to \cite{GritzmannSturmfels} for a detailed analysis.
\end{remark}

\subsection{Sample code}
We conclude with few examples  which illustrate the object  \texttt{HyperplaneArrangement} and its property. Example \ref{01vectors} reports the comparison between the running times 
of the new  algorithm implemented in \polymake  and the brute force algorithm to compute cell decompositions.

\begin{example} The following examples compute the $4! = 24$ cells in the Coxeter arrangement of type A3. The $6$ linear hyperplanes in the arrangements are 
\[
x_i - x_j = 0, \ \ 1 \leq i < j \leq 4.
\] 
\begin{lstlisting}
fan > $A3 = new HyperplaneArrangement(HYPERPLANES=>root_system("A3")->VECTORS->minor(All,~[0]));

fan > $CDA3 = $A3->CELL_DECOMPOSITION;

fan > print $CDA3->N_MAXIMAL_CONES;
24
\end{lstlisting}
Now we compute the $36$ cells in the Linial arrangement \cite{PostnikovStanley} given by the $6$ affine hyperplanes 
\[
x_i - x_j = 1, \ \ 1 \leq i < j \leq 4.
\] 
As explained in Section 2.1, the support cone allows us to deal with affine hyperplanes. We transform the hyperplanes $[a,b] \in \RR^4 \times \RR$ in the projective arrangement $H_{\text{proj}}$ with hyperplanes $[-b,a]$ and then we intersect the latter with the support cone $\supportcone_{H_{\text{proj}}}\ :=\ \{[x_0,x_1,\ldots,x_{5}]\in\RR^{5}\ |\ x_0\ge 0\}$. 
\begin{lstlisting}
fan > $Hyps = new Matrix([[-1,1,-1,0,0],[-1,1,0,-1,0], 
[-1,1,0,0,-1],[-1,0,1,-1,0],[-1,0,1,0,-1],[-1,0,0,1,-1]]);

fan > $Lin = new HyperplaneArrangement(HYPERPLANES=>$Hyps,    
 "SUPPORT.INEQUALITIES"=>[[1,0,0,0,0]]);

fan > $CDLin = $Lin->CELL_DECOMPOSITION;

fan > print $CDLin->N_MAXIMAL_CONES;
36
\end{lstlisting}

\begin{figure}[ht]
\centering
\includegraphics[width=5cm]{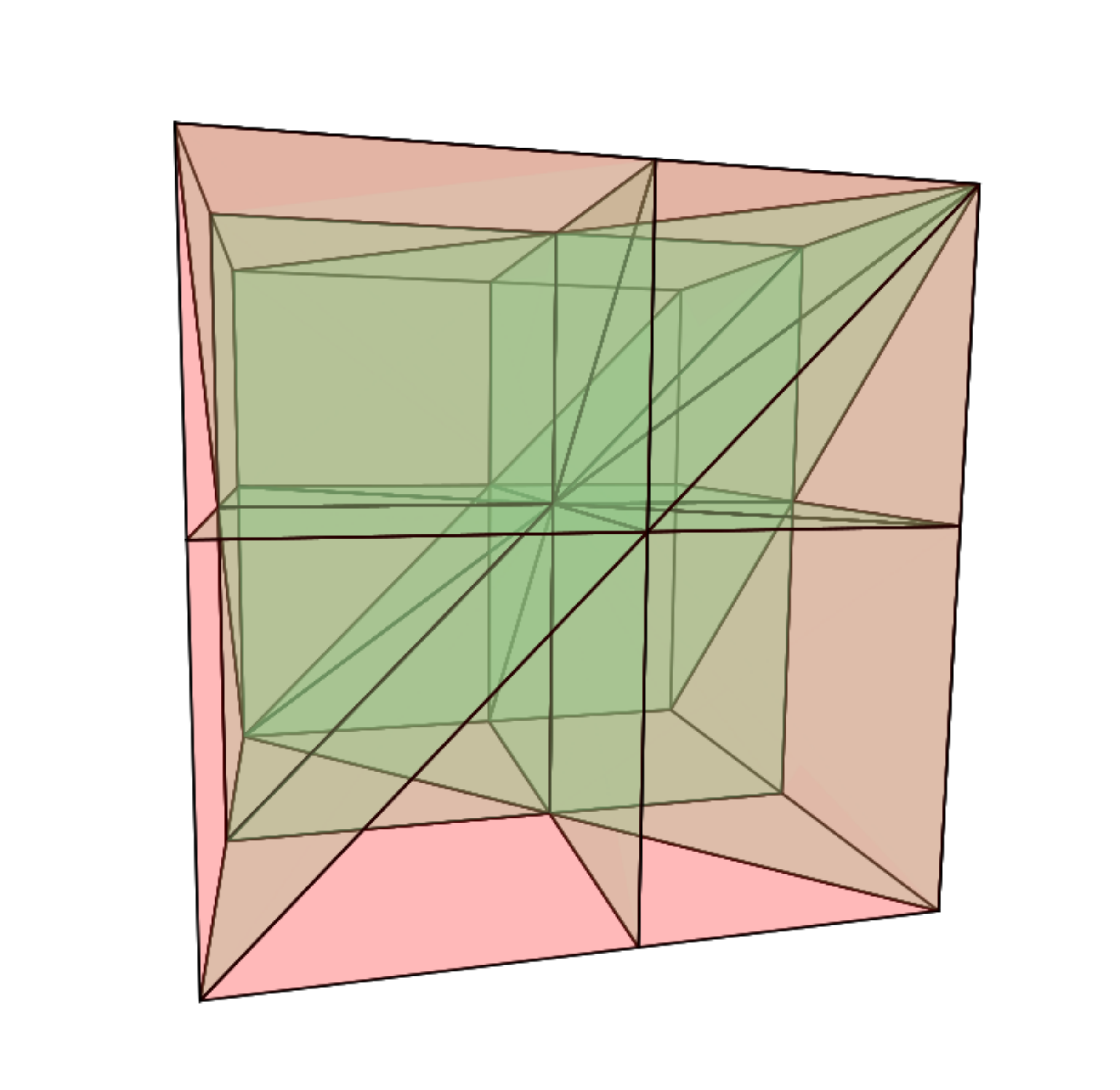}
\caption{The arrangement of type A3.}
\label{fig:A3}
\end{figure}

\end{example}

\begin{example} \label{delPezzo} 
This example is based on \cite{suess19}. Let $X$ be a del Pezzo surface of degree $5$ and $[K_X]$ the class of the canonical divisor. The cone of effective divisors $\overline{\textrm{Eff}}(X)$ is spanned by ten exceptional curves $[C_{ij}]$ indexed by $0 \leq i<j \leq 4$ and  characterized by $[C_{ij}]^2=-1$ and $[C_{ij}]\cdot [K_Y]=-1$. Applying the change of basis $[C_{ij}] =b_i + b_j$, described in \cite[Section 3]{suess19}, we see that the  polytope $P$ given by points in  $\overline{\textrm{Eff}}(X)$ intersecting $[K_X]$  with multiplicity $-1$ coincides with the hypersimplex $\Delta(2,5)$.  In the aforementioned article the author  considers  the cell decomposition of $P$ induced by the hyperplane arrangement defined by 
\[\{ [D] \in \overline{\textrm{Eff}}(X) \, | \, [D] \cdot [C_{ij}] =0\}.\]
The decomposition is used to study the toric topology of the Grassmannian of planes in complex $5$-dimensional space. 

The following code allows one to compute the cell decomposition in \polymake. We first compute  the new pairing in the new basis $b_0, b_1, \dots, b_4$.

\begin{lstlisting}
polytope > $pairing = new Matrix(1/4*ones_matrix(5,5));
polytope > $pairing->row(4) *= -1;
polytope > $pairing->col(4) *= -1;
polytope > for(my $i=0; $i<4; $i++){ $pairing->elem($i,$i) = -3/4; }
\end{lstlisting}

We then introduce the support cone given by the hypersimplex $\Delta(2,5)$ in the new basis
\begin{lstlisting}
polytope > $R = hypersimplex(2,5)->VERTICES;
polytope > $Z = zero_vector(10);
polytope > $M = hypersimplex(2,5)->VERTICES->minor(All,~[0]);
polytope > $M = $M * $pairing;
polytope > $H = $Z|$M;
\end{lstlisting}

Finally, we can compute the cell decomposition. 

\begin{lstlisting}
fan > $HA = new HyperplaneArrangement(HYPERPLANES=>$H, 
"SUPPORT.INPUT_RAYS"=>$R);

fan > print $HA->CELL_DECOMPOSITION->N_RAYS;
15

fan > print $HA->CELL_DECOMPOSITION->N_MAXIMAL_CONES;
27
\end{lstlisting}

\end{example}

\begin{example} \label{01vectors} Let $H$ be the hyperplane arrangement in $\RR^d$ given by the the $2^d-1$ hyperplanes normal to $0\slash 1$-vectors. The number of maximal cones 
in $\Sigma_H$ are known up to $d=8$, see entry A034997 in the Online Encyclopedia of Integer Sequences. We run \polymake implementations of the BFS algorithm described above and the 
brute force alternative. Our results are reported in Table \ref{tab:summary}, where we can see that the BFS algorithm performs better than  already for small values of $d$. 

\begin{table}[h]\centering
\caption{Results and runtimes  for arrangements in Example \ref{01vectors}}
\label{tab:summary}
\setlength{\tabcolsep}{.4em}
\begin{tabular}{crrrrr}
\toprule
$d$ & $\#$ hyperplanes & $\#$ rays & $\#$ maximal cones & time BFS (s) & time brute force (s)\\
\midrule
2 &  3 &    6 &     6 &   0.1 &   0.1\\
3 &  7 &   18 &    32 &   0.6 &   1.2\\
4 & 15 &   90 &   370 &   8.0 & 324.1\\
5 & 31 & 1250 & 11292 & 407.6 &    --\\
\bottomrule
\end{tabular}
\end{table}

\end{example}

\printbibliography

\end{document}